# Real Hypersurfaces of Type A in Complex Two-Plane Grassmannians Related to The Reeb Vector Field

Ruenn-Huah Lee[a] and Tee-How Loo[b]

[a,b]Institute of Mathematical Sciences, Faculty of Science
University of Malaya, 50603 Kuala Lumpur, Malaysia
*rhlee063@siswa.um.edu.my, looth@um.edu.my*

**Abstract.** Y. J. Suh and H. Lee (Bull. Korean. Math. Soc. **47**, 551-561 (2010)) characterized real hypersurfaces $M$ of type $B$ by the invariance of vector bundle $JTM^\perp$ under the shape operator and the orthogonality of $JTM^\perp$ and $\mathrm{J}\,TM^\perp$, where $TM^\perp$, $J$ and $\mathrm{J}$ are the normal bundle of $M$, Kähler structure and Quaternionic Kähler structure of $G_2(C^{m+2})$ respectively. In this paper, we characterize real hypersurfaces $M$ of type $A$ by the invariance of the vector bundle $JTM^\perp$ under the shape operator with the Reeb vector field in $\mathrm{J}\,TM^\perp$.

## INTRODUCTION

The beauty of real hypersurfaces of complex two-plane Grassmannians was not discovered and fully appreciated until J. Berndt provided a thorough study on the Riemannian geometry of the complex two-plane Grassmannians [1]. The complex two-plane grassmannian has some remarkable properties and structures. The most notable one being the fact that it is the unique compact irreducible Riemannian symmetric space with both a Kähler structure $J$ and a quaternionic Kähler structure $\mathrm{J}$. These geometric structures induce a local almost contact 3-structure $(\phi_a, \xi_a, \eta_a)$, $a \in \{1,2,3\}$ as well as an almost contact structure $(\phi, \xi, \eta)$ on its real hypersurfaces $M$.

The development of the theory of real hypersurfaces in complex two-plane Grassmannians has been at its prime after the characterization of real hypersurfaces of type $A$ and type $B$ given by J. Berndt and Y. J. Suh in the late 90s [2]. In the following, we denote by $\mathrm{D}^\perp$, the space spanned by $\{\xi_1, \xi_2, \xi_3\}$:

**Theorem 1.1.** ([2]) Let $M$ be a connected real hypersurface in $G_2(C^{m+2})$, $m \geq 3$. Then both $\xi$ and $\mathrm{D}^\perp$ are invariant under the shape operator of $M$ if and only if
   (A) $M$ is an open part of a tube around a totally geodesic $G_2(C^{m+1})$ of $G_2(C^{m+2})$, or
   (B) $m$ is even, say $m = 2n$, and $M$ is an open part of a tube around a totally geodesic $\mathrm{H}P^n$ in $G_2(C^{m+2})$.

We say that a real hypersurface $M$ in $G_2(C^{m+2})$ is of *type A* if it satisfies the first property in the characterization theorem given above. On the other hand, $M$ is said to be of type $B$ if it satisfies all properties in part (B). A connected orientable real hypersurface $M$ in $G_2(C^{m+2})$ is said to be *Hopf* if the Reeb vector field $\xi$ is invariant under the shape operator of $M$. The following theorems were also given in the same paper, provide sufficient conditions of being a real hypersurface of type $A$ and type $B$.

**Theorem 1.2.** ([2]) Let $M$ be a connected real hypersurface of $G_2(C^{m+2})$, $m \geq 3$. If $A\mathrm{D} \subset \mathrm{D}$, $A\xi = \alpha\xi$ and $\xi$ is tangent to $\mathrm{D}^\perp$, then $M$ is an open part of a tube around a totally geodesic $G_2(C^{m+1})$ in $G_2(C^{m+2})$.

**Theorem 1.3**. ([2]) Let $M$ be a connected real hypersurface of $G_2(C^{m+2})$, $m \geq 3$. If $A\mathrm{D} \subset \mathrm{D}$, $A\xi = \alpha\xi$ and $\xi$ is tangent to $\mathrm{D}$, then $m = 2n$ and $M$ is an open part of a tube around a totally geodesic $\mathrm{H}P^n$ in $G_2(C^{m+2})$.

These characterizations of real hypersurfaces of type $A$ and type $B$ look extremely similar. Recently, H. Lee and Y. J. Suh shown that the invariance of $\mathfrak{D}$ under $A$ is in fact a consequence of being a Hopf hypersurface with the Reeb vector field belongs to $\mathfrak{D}$ [3].

**Theorem 1.4.** ([3]) Let $M$ be a connected orientable Hopf hypersurface in $G_2(\mathbb{C}^{m+2})$, $m \geq 3$. Then the Reeb vector $\xi$ belongs to the distribution $\mathfrak{D}$ if and only if $M$ is locally congruent to an open part of a real hypersurface of type $B$.

The preceding theorem improves Theorem 1.3 significantly. It is reasonable to wonder if it is possible to weaken the hypothesis in Theorem 1.2 in a similar manner. In this paper, we give a new characterization of real hypersurfaces of type $A$ as motivated by the theorem given above, but instead of assuming that $M$ is Hopf, we assume that the distribution $\mathfrak{D}$ is invariant under the shape operator $A$. We proved the following theorem:

**Theorem 1.5.** Let $M$ be a connected orientable real hypersurface in $G_2(\mathbb{C}^{m+2})$, $m \geq 3$. If $A\mathfrak{D} \subset \mathfrak{D}$ and $\xi \in \mathfrak{D}^\perp$, then $M$ is Hopf.

By Theorem 1.2 and Theorem 1.5, we have:

**Theorem 1.6.** Let $M$ be a connected orientable real hypersurface in $G_2(\mathbb{C}^{m+2})$, $m \geq 3$. If $A\mathfrak{D} \subset \mathfrak{D}$ and $\xi \in \mathfrak{D}^\perp$, then $M$ is an open part of a tube around a totally geodesic $G_2(\mathbb{C}^{m+1})$ in $G_2(\mathbb{C}^{m+2})$.

## REAL HYPERSURFACES IN $G_2(\mathbb{C}^{m+2})$

In this section, we summarize and list out some important formulae as well as well-known results in the theory of real hypersurfaces in complex two-plane Grassmannians. We begin this section with the notations that we are going to use throughout this paper and also some basic properties of complex two-plane Grassmannians.

First, we denote the set of all complex 2-dimensional linear subspaces of $\mathbb{C}^{m+2}$ by $G_2(\mathbb{C}^{m+2})$. The Riemannian metric is denoted by $g$. Next, we denote by $J$ and $\mathfrak{J}$ the Kähler structure and the quaternionic Kähler structure on $G_2(\mathbb{C}^{m+2})$ respectively. For each $x \in G_2(\mathbb{C}^{m+2})$, we denote by $\{J_1, J_2, J_3\}$ a canonical local basis of $\mathfrak{J}$ on a neighborhood U of $x$ in $G_2(\mathbb{C}^{m+2})$, that is, each $J_a$ is a local almost Hermitian structure such that

$$J_a J_{a+1} = J_{a+2} = -J_{a+1} J_a, \qquad a \in \{1, 2, 3\}.$$

The indices in the preceding equation is taken modulo three. The Levi-Civita connection of $G_2(\mathbb{C}^{m+2})$ is denoted by $\bar{\nabla}$. Since $\mathfrak{J}$ is parallel with respect to the Levi-Civita connection $\bar{\nabla}$ of $G_2(\mathbb{C}^{m+2})$, for any canonical local basis $\{J_1, J_2, J_3\}$ of $\mathfrak{J}$, there exist three local 1-forms $q_1, q_2, q_3$ such that

$$\bar{\nabla}_X J_a = q_{a+2}(X) J_{a+1} - q_{a+1}(X) J_{a+2},$$

for each vecor field $X$ on $G_2(\mathbb{C}^{m+2})$, where $\bar{\nabla}$ is the Levi-Civita connection on $G_2(\mathbb{C}^{m+2})$. The Kähler structure $J$ and the quaternionic Kähler structure $\mathfrak{J}$ are related by

$$J J_a = J_a J \qquad \text{and} \qquad \text{Trace}(J J_a) = 0, \tag{1}$$

for all $a \in \{1, 2, 3\}$. Note that all the indices are taken modulo three. Next, the Riemannian curvature tensor $\bar{R}$ of $G_2(\mathbb{C}^{m+2})$ is locally given by

$$\begin{aligned}
\bar{R}(X,Y)Z = {} & g(Y,Z)X - g(X,Z)Y \\
& + g(JY,Z)JX - g(JX,Z)JY - 2g(JX,Y)JZ \\
& + \sum_{a=1}^{3}\{g(J_aY,Z)J_aX - g(J_aX,Z)J_aY - 2g(J_aX,Y)J_aZ \\
& + g(JJ_aY,Z)JJ_aX - g(JJ_aX,Z)JJ_aY\},
\end{aligned} \qquad (2)$$

for all $X,Y,Z \in T_x(G_2(C^{m+2}))$.

For a nonzero vector $X \in T_x(G_2(C^{m+2}))$, we denote by $CX = \text{span}\{X, JX\}$, $\mathsf{J}X = \{J'X \mid J' \in \mathsf{J}_x\}$, $HX = \text{span}\{X\} \oplus \mathsf{J}X$, and $HCX$ the subspace spanned by $HX$ and $HJX$. If $JX \in \mathsf{J}X$, we denote by $C^\perp X$ the orthogonal complement of $CX$ in $HX$.

Let $M$ be a connected, oriented real hypersurface isometrically immersed in $G_2(C^{m+2})$, $m \geq 3$, and $N$ a unit normal vector field on $M$. The Riemannian metric on $M$ will also be denoted by $g$ as there will be no confusion occured. A canonical local basis $\{J_1, J_2, J_3\}$ of $\mathsf{J}$ on $G_2(C^{m+2})$ induces a local almost contact metric 3-structure $(\phi_a, \xi_a, \eta_a, g)$ on $M$ by

$$J_a X = \phi_a X + \eta_a(X) N, \qquad J_a N = -\xi_a, \qquad \eta_a(X) = g(X, \xi_a),$$

for all vector field $X$ tangent to $M$. It follows that

$$\phi_a \phi_{a+1} - \xi_a \otimes \eta_{a+1} = \phi_{a+2} = -\phi_{a+1}\phi_a + \xi_{a+1} \otimes \eta_a,$$
$$\phi_a \xi_{a+1} = \xi_{a+2} = -\phi_{a+1}\xi_a.$$

Let $(\phi, \xi, \eta, g)$ be the almost contact metric structure on $M$ induced by $J$, that is,

$$JX = \phi X + \eta(X)N, \qquad JN = -\xi, \qquad \eta(X) = g(X, \xi).$$

A real hypersurface $M$ is said to be *Hopf* if the Reeb vector field $\xi$ is principal.

It follows from (1) that the two structures $(\phi, \xi, \eta, g)$ and $(\phi_a, \xi_a, \eta_a, g)$ are related as follows

$$\phi_a\phi - \xi_a \otimes \eta = \phi\phi_a - \xi \otimes \eta_a; \qquad \phi\xi_a = \phi_a\xi.$$

Next, we denote by $\nabla$ the Levi-Civita connection and $A$ the shape operator on $M$. Then

$$\begin{aligned}
(\nabla_X \phi)Y &= \eta(Y)AX - g(AX,Y)\xi, \qquad \nabla_X\xi = \phi AX, \\
(\nabla_X \phi_a)Y &= \eta_a(Y)AX - g(AX,Y)\xi_a + q_{a+2}(X)\phi_{a+1}Y - q_{a+1}(X)\phi_{a+2}Y, \\
\nabla_X \xi_a &= \phi_a AX + q_{a+2}(X)\xi_{a+1} - q_{a+1}(X)\xi_{a+2},
\end{aligned}$$

for any $X, Y$ tangent to $M$.

Finally we state some well-known results.

**Lemma 2.1.** ([4]) Let $M$ be a real hypersurface in $G_2(C^{m+2})$, $m \geq 3$. If $\xi$ is tangent to $\mathsf{D}$, then $A\phi\xi_a = 0$, for $a \in \{1,2,3\}$.

**Theorem 2.2.** ([2]) Let $M$ be a real hypersurface of type $A$ in $G_2(C^{m+2})$, $m \geq 3$. Then $\xi \in \mathsf{D}^\perp$ at each point of $M$. Suppose $J_1 \in \mathsf{J}$ such that $J_1 N = JN$. Then $M$ has three (if $r = \pi/2\sqrt{8}$) or four (otherwise) distinct constant principal curvatures

$$\alpha = \sqrt{8}\cot(\sqrt{8}r), \qquad \beta = \sqrt{2}\cot(\sqrt{2}r), \qquad \lambda = -\sqrt{2}\tan(\sqrt{2}r), \qquad \mu = 0$$

with some $r \in ]0, \pi/\sqrt{8}[$. The corresponding multiplicities are

$$m(\alpha) = 1, \qquad m(\beta) = 2, \qquad m(\lambda) = 2m - 2 = m(\mu)$$

and the corresponding eigenspaces are

$$\begin{aligned}T_\alpha &= span\{\xi\} \\ T_\beta &= C^\perp \xi \\ T_\lambda &= \{X : X \perp H\xi, JX = J_1 X\} \\ T_\mu &= \{X : X \perp H\xi, JX = -J_1 X\}\end{aligned}$$

Next, we give a short treatment on a symmetric tensor field $\theta_a$ introduced in [4], which is essential in the proof of our main theorem. Let $M$ be a real hypersurface in $G_2(C^{m+2})$, $m \geq 3$. Corresponding to each canonical local basis $\{J_1, J_2, J_3\}$ of J , we define a local endomorphism $\theta_a$ on $TM$ by

$$\theta_a X := \tan(JJ_a) = \phi_a \phi X - \eta(X)\xi_a = \phi\phi_a X - \eta_a(X)\xi.$$

It follows from (2) that the equation of Gauss is given by

$$\begin{aligned}R(X,Y)Z = &\ g(Y,Z)X - g(X,Z)Y + g(AY,Z)AX - g(AX,Z)AY \\ &+ g(\phi Y,Z)\phi X - g(\phi X,Z)\phi Y - 2g(\phi X,Y)\phi Z \\ &+ \sum_{a=1}^{3}\{g(\phi_a Y,Z)\phi_a X - g(\phi_a X,Z)\phi_a Y - 2g(\phi_a X,Y)\phi_a Z \\ &+ g(\theta_a Y,Z)\theta_a X - g(\theta_a X,Z)\theta_a Y\}.\end{aligned}$$

Similarly, by adopting the symmetric tensor field $\theta_a$, the Codazzi equation is given by

$$\begin{aligned}(\nabla_X A)Y - (\nabla_Y A)X = &\ \eta(X)\phi Y - \eta(Y)\phi X - 2g(\phi X,Y)\xi \\ &+ \sum_{a=1}^{3}(\eta_a(X)\phi_a Y - \eta_a(Y)\phi_a X - 2g(\phi_a X,Y)\xi_a \\ &+ \eta_a(\phi X)\theta_a Y - \eta_a(\phi Y)\theta_a X).\end{aligned} \qquad (3)$$

For each $x \in M$, we define a subspace $H^\perp$ of $T_x M$ by

$$H^\perp := span\{\xi, \xi_1, \xi_2, \xi_3, \phi\xi_1, \phi\xi_2, \phi\xi_3\}.$$

Let H be the orthogonal complement of $HC\xi$ in $T_x(G_2(C^{m+2}))$. Then, $T_x M = H \oplus H^\perp$ and H is invariant under $\phi, \phi_a$ and $\theta_a$. Moreover, $\theta_{a|H}$ has two eigenvalues: 1 and $-1$. The following lemmas summarize some important properties and facts about $\theta_a$:

**Lemma 2.4.** ([4])
  (a) $\theta_a$ is symmetric,
  (b) Trace$(\theta_a) = \eta(\xi_a)$,

(c) $\theta_a^2 X = X - g(X, \phi\xi_a)\phi\xi_a$, for all $X \in TM$,

(d) $\theta_a \xi = -\xi_a$; $\quad \theta_a \xi_a = -\xi$; $\quad \theta_a \phi\xi_a = \eta(\xi_a)\phi\xi_a$,

(e) $\theta_a \xi_{a+1} = \phi\xi_{a+2} = -\theta_{a+1}\xi_a$,

(f) $\theta_a \phi\xi_{a+1} = -\xi_{a+2} + \eta(\xi_{a+1})\phi\xi_a$,

(g) $\theta_{a+1}\phi\xi_a = \xi_{a+2} + \eta(\xi_a)\phi\xi_{a+1}$.

**Lemma 2.5.** ([4]) Let $H_a(\varepsilon)$ be the eigenspace corresponding to eigenvalue $\varepsilon$ of $\theta_{a|H}$. Then

(a) $\theta_{a|H}$ has two eigenvalues $\varepsilon = \pm 1$,

(b) $\phi H_a(\varepsilon) = H_a(\varepsilon)$,

(c) $\theta_b H_a(\varepsilon) = H_a(-\varepsilon)$, for $a \neq b$,

(d) $\dim H_a(1) = \dim H_a(-1)$ is even,

(e) $\phi_b H_a(\varepsilon) = H_a(-\varepsilon)$, for $a \neq b$.

Note that $\xi \in D^\perp$ if and only if $\dim H^\perp = 3$, that is, when $\xi \in D^\perp$, we have $H^\perp = D^\perp$. One may refer to [4] for more details.

# PROOF OF THEOREM 1.5

The assumption that $AD \subset D$ also implies $AD^\perp \subset D^\perp$. Hence, by a suitable choice of canonical local basis $\{J_1, J_2, J_3\}$ for J, the vector fields $\xi_1, \xi_2, \xi_3$ are principal everywhere, say $A\xi_a = \beta_a \xi_a$, for $a = 1, 2, 3$.

From the Codazzi equation (3), we have:

$$g((\nabla_X A)Y - (\nabla_Y A)X, \xi_a) = -2\eta(\xi_a)g(\phi X, Y) - 2g(\phi_a X, Y)$$
$$+2\eta(X)\eta_a(\phi Y) - 2\eta(Y)\eta_a(\phi X)$$
$$+2\eta_{a+1}(X)\eta_{a+2}(Y) - 2\eta_{a+1}(Y)\eta_{a+2}(X)$$
$$+2\eta_{a+1}(\phi X)\eta_{a+2}(\phi Y) - 2\eta_{a+1}(\phi Y)\eta_{a+2}(\phi X),$$

for all vector fields $X, Y$ tangent to $M$. On the other hand, since the shape operator $A$ is symmetric, $\nabla_X A$ is also symmetric, that is, $g((\nabla_X A)Y, Z) = g(Y, (\nabla_X A)Z)$. Due to the fact that the canonical local basis is chosen in a way that $\xi_a$, $a = 1, 2, 3$ are principal, we have

$$g((\nabla_X A)Y - (\nabla_Y A)X, \xi_a) = g((\nabla_X A)\xi_a, Y) - g((\nabla_Y A)\xi_a, X)$$
$$= X(\beta_a)\eta_a(Y) - Y(\beta_a)\eta_a(X)$$
$$+ \beta_a g((\phi_a A + A\phi_a)X, Y) - 2g(A\phi_a AX, Y)$$
$$+ (\beta_a - \beta_{a+1})[q_{a+2}(X)\eta_{a+1}(Y) - q_{a+2}(Y)\eta_{a+1}(X)]$$
$$- (\beta_a - \beta_{a+2})[q_{a+1}(X)\eta_{a+2}(Y) - q_{a+1}(Y)\eta_{a+2}(X)].$$

Combining these equations, we have

$$-2\eta(\xi_a)g(\phi X, Y) - 2g(\phi_a X, Y) + 2\eta(X)\eta_a(\phi Y) - 2\eta(Y)\eta_a(\phi X)$$
$$+2\eta_{a+1}(X)\eta_{a+2}(Y) - 2\eta_{a+1}(Y)\eta_{a+2}(X) + 2\eta_{a+1}(\phi X)\eta_{a+2}(\phi Y) - 2\eta_{a+1}(\phi Y)\eta_{a+2}(\phi X)$$
$$= X(\beta_a)\eta_a(Y) - Y(\beta_a)\eta_a(X) + \beta_a g((\phi_a A + A\phi_a)X, Y) - 2g(A\phi_a AX, Y) \qquad (4)$$
$$+(\beta_a - \beta_{a+1})[q_{a+2}(X)\eta_{a+1}(Y) - q_{a+2}(Y)\eta_{a+1}(X)]$$
$$-(\beta_a - \beta_{a+2})[q_{a+1}(X)\eta_{a+2}(Y) - q_{a+1}(Y)\eta_{a+2}(X)].$$

By taking $Y = \xi_a$, we obtain

$$X(\beta_a) = \xi_a(\beta_a)\eta_a(X) - 4\eta(\xi_a)\eta_a(\phi X) + 2\eta(\xi_{a+1})\eta_{a+1}(\phi X) + 2\eta(\xi_{a+2})\eta_{a+2}(\phi X)$$
$$+ (\beta_a - \beta_{a+1})q_{a+2}(\xi_a)\eta_{a+1}(X) - (\beta_a - \beta_{a+2})q_{a+1}(\xi_a)\eta_{a+2}(X).$$

$Y(\beta_a)$ could be obtained in a similar manner. Some of these computations can be found in [2]. The following lemma is obtained by substituting both $X(\beta_a)$ and $Y(\beta_a)$ back into (4).

**Lemma 3.1.** Let $M$ be a connected orientable hypersurface in $G_2(\mathbb{C}^{m+2})$, $m \geq 3$. Suppose $A\xi_a = \beta_a \xi_a$, for $a = 1, 2, 3$, then

$$2\eta(\xi_a)g(\phi X, Y) + 2g(\phi_a X, Y) + \beta_a g((\phi_a A + A\phi_a)X, Y) - 2g(A\phi_a AX, Y)$$
$$= +2\eta(X)\eta_a(\phi Y) - 2\eta(Y)\eta_a(\phi X)$$
$$+ 2\eta_{a+1}(X)\eta_{a+2}(Y) - 2\eta_{a+1}(Y)\eta_{a+2}(X) + 2\eta_{a+1}(\phi X)\eta_{a+2}(\phi Y) - 2\eta_{a+1}(\phi Y)\eta_{a+2}(\phi X)$$
$$+ 2\eta_a(Y)[2\eta(\xi_a)\eta_a(\phi X) - \eta(\xi_{a+1})\eta_{a+1}(\phi X) - \eta(\xi_{a+2})\eta_{a+2}(\phi X)]$$
$$- 2\eta_a(X)[2\eta(\xi_a)\eta_a(\phi Y) - \eta(\xi_{a+1})\eta_{a+1}(\phi Y) - \eta(\xi_{a+2})\eta_{a+2}(\phi Y)] \quad (5)$$
$$+ (\beta_a - \beta_{a+1})q_{a+2}(\xi_a)[\eta_a(X)\eta_{a+1}(Y) - \eta_{a+1}(X)\eta_a(Y)]$$
$$+ (\beta_a - \beta_{a+2})q_{a+1}(\xi_a)[\eta_{a+2}(X)\eta_a(Y) - \eta_a(X)\eta_{a+2}(Y)]$$
$$- (\beta_a - \beta_{a+1})[q_{a+2}(X)\eta_{a+1}(Y) - q_{a+2}(Y)\eta_{a+1}(X)]$$
$$+ (\beta_a - \beta_{a+2})[q_{a+1}(X)\eta_{a+2}(Y) - q_{a+1}(Y)\eta_{a+2}(X)],$$

for all $X, Y$ tangent to $M$.

We further assume that $\xi \in \mathsf{D}^\perp$. For any $X \in \mathsf{D}$, since $\mathsf{D}$ is invariant under $A, \phi$ and $\phi_a$, from equation (5), we have

$$0 = 2\eta(\xi_a)g(\phi X, Y) + 2g(\phi_a X, Y) + \beta_a g((\phi_a A + A\phi_a)X, Y) - 2g(A\phi_a AX, Y)$$
$$- (\beta_a - \beta_{a+1})q_{a+2}(X)\eta_{a+1}(Y) + (\beta_a - \beta_{a+2})q_{a+1}(X)\eta_{a+2}(Y), \quad (6)$$

for all $Y$ tangent to $M$ and $a = 1, 2, 3$. Hence, we have

**Lemma 3.2.** Let $M$ be a connected orientable hypersurface in $G_2(\mathbb{C}^{m+2})$, $m \geq 3$. Suppose $A\mathsf{D}^\perp \subset \mathsf{D}^\perp$ and $\xi \in \mathsf{D}^\perp$, if $A\xi_a = \beta_a \xi_a$, for $a = 1, 2, 3$, then

$$2\eta(\xi_a)\phi X + 2\phi_a X + \beta_a(\phi_a A + A\phi_a)X - 2A\phi_a AX = 0, \quad (7)$$

for all $X \in \mathsf{D}$.

**Lemma 3.3.** Let $M$ be a connected orientable hypersurface in $G_2(\mathbb{C}^{m+2})$, $m \geq 3$. Suppose $A\mathsf{D}^\perp \subset \mathsf{D}^\perp$ and $\xi \in \mathsf{D}^\perp$, if $A\xi_a = \beta_a \xi_a$ for $a = 1, 2, 3$, then for all $X \in \mathsf{D}$, we have $\phi_a A\phi_a AX = A\phi_a A\phi_a X$.

*Proof.* Note that if $X \in \mathsf{D}$, then $\phi_a X \in \mathsf{D}$. Next, applying $\phi_a$ on both sides of (7) and replacing $X$ by $\phi_a X$ in (7) give

$$2\eta(\xi_a)\theta_a X - 2X - \beta_a AX + \beta_a \phi_a A\phi_a X - 2\phi_a A\phi_a AX = 0$$

and

$$2\eta(\xi_a)\theta_a X - 2X + \beta_a\phi_a A\phi_a X - \beta_a AX - 2A\phi_a A\phi_a X = 0$$

respectively. Hence, $\phi_a A\phi_a AX = A\phi_a A\phi_a X$, for all $X \in \mathsf{D}$. □

With the assumption that $\xi \in \mathsf{D}^\perp$, we have $\mathsf{D} = \mathsf{H}$ and $A\mathsf{H} \subset \mathsf{H}$. By Lemma 3.3, there exist common orthonormal eigenvectors $X_1, \cdots, X_{4m-4} \in \mathsf{H}$ of $A$ and $\phi_1 A\phi_1$. It follows that

$$AX_j = \lambda_j X_j, \quad \text{and} \quad A\phi_1 X_j = \mu_j \phi_1 X_j. \tag{8}$$

Using these in (7), we have

$$2\eta(\xi_1)\phi X_j + (2 + \lambda_j \beta_1 + \mu_j \beta_1 - 2\lambda_j \mu_j)\phi_1 X_j = 0.$$

Since $\xi \in \mathsf{D}^\perp$, we may suppose $\eta(\xi_1) \neq 0$. Then

$$0 = \phi X_j + \varepsilon \phi_1 X_j, \tag{9}$$

where $\varepsilon = \dfrac{2 + \beta_1(\lambda_j + \mu_j) - 2\lambda_j \mu_j}{2\eta(\xi_1)}$. Taking $\phi_1$ on both sides of (9) gives

$$0 = \theta_1 X_j - \varepsilon X_j.$$

By Lemma 2.5 (a), $\varepsilon \in \{1, -1\}$. Without loss of generality, we can assume that

$$X_1, \cdots, X_{2m-2} \in \mathsf{H}_1(1) \quad \text{and} \quad X_{2m-1}, \cdots, X_{4m-4} \in \mathsf{H}_1(-1),$$

where $\mathsf{H}_1(1)$ and $\mathsf{H}_1(-1)$ are the eigenspace of $\theta_1|_\mathsf{H}$ corresponding to the eigenvalue 1 and $-1$ respectively. Consequently, (8) implies

$$A\mathsf{H}_1(1) \subset \mathsf{H}_1(1)$$

and hence

$$\phi_2 A\phi_2 \mathsf{H}_1(1) \subset \mathsf{H}_1(1).$$

Thus, if we take $a = 2$ in Lemma 3.3, then there exists orthonormal vectors $\bar{X}_1, \cdots, \bar{X}_{2m-2} \in \mathsf{H}_1(1)$ such that

$$A\bar{X}_j = \bar{\lambda}_j \bar{X}_j \quad \text{and} \quad A\phi_2 \bar{X}_j = \bar{\mu}_j \phi_2 \bar{X}_j.$$

From (7), we have

$$2\eta(\xi_2)\phi\bar{X}_j + (2 + \bar{\lambda}_j \beta_2 + \bar{\mu}_j \beta_2 - 2\bar{\lambda}_j \bar{\mu}_j)\phi_2 \bar{X}_j = 0.$$

Since $\bar{X}_j \in \mathsf{H}_1(1)$, $\phi\bar{X}_j \in \mathsf{H}_1(1)$ and $\phi_2 \bar{X}_j \in \mathsf{H}_1(-1)$, we have

$$\eta(\xi_2) = 0.$$

In a similar manner, we obtain $\eta(\xi_3) = 0$. Thus, we have $\eta(\xi_1) = \pm 1$ or $\xi = \pm\xi_1$. As a result, we have shown that $A\xi = \alpha\xi$ with $\alpha = \beta_1$. This completes the proof of Theorem 1.5.